\documentclass[10pt]{amsart}

\usepackage{amsmath,amsthm,amsfonts,amscd,amssymb}
\usepackage{hyperref,wasysym}
\usepackage{verbatim}
\usepackage{graphicx}
\usepackage{tikz}

\newtheorem{thm}{Theorem}[section]
\newtheorem{cor}[thm]{Corollary}
\newtheorem{lem}[thm]{Lemma}
\newtheorem{prop}[thm]{Proposition}

\theoremstyle{definition}

\theoremstyle{remark}
\newtheorem{rem}{Remark}[section]

\begin{document}

\title{On lattice extensions}

\author{Maxwell Forst}
\author{Lenny Fukshansky}\thanks{Fukshansky was partially supported by the Simons Foundation grant \#519058}

\address{Institute of Mathematical Sciences, Claremont Graduate University, Claremont, CA 91711}
\email{maxwell.forst@cgu.edu}
\address{Department of Mathematics, 850 Columbia Avenue, Claremont McKenna College, Claremont, CA 91711}
\email{lenny@cmc.edu}


\begin{abstract} A lattice $\Lambda$ is said to be an extension of a sublattice $L$ of smaller rank if $L$ is equal to the intersection of $\Lambda$ with the subspace spanned by $L$. The goal of this paper is to initiate a systematic study of the geometry of lattice extensions. We start by proving the existence of a small-determinant extension of a given lattice, and then look at successive minima and covering radius. To this end, we investigate extensions (within an ambient lattice) preserving the successive minima of the given lattice, as well as extensions preserving the covering radius. We also exhibit some interesting arithmetic properties of deep holes of planar lattices.
\end{abstract}

\maketitle

\def\A{{\mathcal A}}
\def\a{{\mathfrak a}}
\def\B{{\mathcal B}}
\def\C{{\mathcal C}}
\def\D{{\mathcal D}}
\def\F{{\mathcal F}}
\def\x{{\mathcal H}}
\def\I{{\mathcal I}}
\def\J{{\mathcal J}}
\def\K{{\mathcal K}}
\def\L{{\mathcal L}}
\def\M{{\mathcal M}}
\def\N{{\mathcal N}}
\def\O{{\mathcal O}}
\def\R{{\mathcal R}}
\def\s{{\mathcal S}}
\def\V{{\mathcal V}}
\def\W{{\mathcal W}}
\def\X{{\mathcal X}}
\def\Y{{\mathcal Y}}
\def\H{{\mathcal H}}
\def\Z{{\mathcal Z}}
\def\OO{{\mathcal O}}
\def\PP{{\mathcal P}}
\def\BB{{\mathbb B}}
\def\cee{{\mathbb C}}
\def\pee{{\mathbb P}}
\def\que{{\mathbb Q}}
\def\real{{\mathbb R}}
\def\zed{{\mathbb Z}}
\def\hyp{{\mathbb H}}
\def\aa{{\mathfrak a}}
\def\HH{{\mathfrak H}}
\def\qbar{{\overline{\mathbb Q}}}
\def\eps{{\varepsilon}}
\def\ahat{{\hat \alpha}}
\def\bhat{{\hat \beta}}
\def\gt{{\tilde \gamma}}
\def\h{{\tfrac12}}
\def\be{{\boldsymbol e}}
\def\bei{{\boldsymbol e_i}}
\def\bff{{\boldsymbol f}}
\def\ba{{\boldsymbol a}}
\def\bb{{\boldsymbol b}}
\def\bc{{\boldsymbol c}}
\def\bm{{\boldsymbol m}}
\def\bk{{\boldsymbol k}}
\def\bi{{\boldsymbol i}}
\def\bl{{\boldsymbol l}}
\def\bq{{\boldsymbol q}}
\def\bu{{\boldsymbol u}}
\def\bt{{\boldsymbol t}}
\def\bs{{\boldsymbol s}}
\def\bv{{\boldsymbol v}}
\def\bw{{\boldsymbol w}}
\def\bx{{\boldsymbol x}}
\def\bX{{\boldsymbol X}}
\def\bz{{\boldsymbol z}}
\def\bwy{{\boldsymbol y}}
\def\bY{{\boldsymbol Y}}
\def\bL{{\boldsymbol L}}
\def\baa{{\boldsymbol\alpha}}
\def\bbb{{\boldsymbol\beta}}
\def\bet{{\boldsymbol\eta}}
\def\bxi{{\boldsymbol\xi}}
\def\bo{{\boldsymbol 0}}
\def\bol{{\boldkey 1}_L}
\def\ep{\varepsilon}
\def\p{\boldsymbol\varphi}
\def\q{\boldsymbol\psi}
\def\rank{\operatorname{rank}}
\def\aut{\operatorname{Aut}}
\def\lcm{\operatorname{lcm}}
\def\sgn{\operatorname{sgn}}
\def\spn{\operatorname{span}}
\def\md{\operatorname{mod}}
\def\Norm{\operatorname{Norm}}
\def\dim{\operatorname{dim}}
\def\det{\operatorname{det}}
\def\Vol{\operatorname{Vol}}
\def\rk{\operatorname{rk}}
\def\Gal{\operatorname{Gal}}
\def\WR{\operatorname{WR}}
\def\WO{\operatorname{WO}}
\def\GL{\operatorname{GL}}
\def\pr{\operatorname{pr}}

\section{Introduction and statement of results}
\label{intro}

Let $n \geq m \geq 2$ be integers and let $\Lambda$ be a lattice of rank $m$ in~$\real^n$, then
$$\Lambda = A\zed^m$$
for an $n \times m$ basis matrix $A$ of rank $m$ and {\it determinant} of $\Lambda$ is
$$\det \Lambda = \sqrt{ \det (A^{\top} A) },$$
which is its co-volume  in $\spn_{\real} \Lambda$: this definition is independent of the choice of the basis matrix $A$ for $\Lambda$. Let $B_m(1)$ be the $m$-dimensional unit ball centered at $\bo$ in $\spn_{\real} \Lambda$ and write $\omega_m$ for the $m$-dimensional volume of $B_m(1)$. Then $B_m(r)$ is a ball of radius $r$ and volume $\omega_m r^m$. We briefly recall the standard notation from the geometry of numbers (see~\cite{lek} for the detailed exposition of the subject). First, the {\it successive minima} of $\Lambda$ are real numbers
$$0 < \lambda_1(\Lambda) \leq \dots \leq \lambda_m(\Lambda),$$
given by
$$\lambda_i(\Lambda) = \min \left\{ r \in \real : \dim_{\real} \spn_{\real} \left( B_m(r) \cap \Lambda \right) \geq i \right\}.$$
Then $\frac{1}{2} \lambda_1(\Lambda)$ is the radius of a ball in the sphere packing associated to $\Lambda$ and the product of successive minima is bounded as follows by the Minkowski's Successive Minima Theorem:
$$\frac{2^m \det \Lambda}{m!\ \omega_m} \leq \prod_{i=1}^m \lambda_i(\Lambda) \leq \frac{2^m \det \Lambda}{\omega_m}.$$
Additionally, the {\it covering radius} of $\Lambda$ is defined as 
$$\mu(\Lambda) = \min \left\{ r \in \real : \Lambda + B_m(r) = \spn_{\real} \Lambda \right\}.$$
The classical inequality of Jarnik asserts that
$$\frac{1}{2} \lambda_m(\Lambda) \leq \mu(\Lambda) \leq \frac{1}{2} \sum_{i=1}^m \lambda_i(\Lambda).$$

Now, let $L \subset \Lambda$ be a sublattice of rank $k < m$. We say that $\Lambda$ is an {\it extension} lattice of $L$ if
$$\Lambda \cap \spn_{\real} L = L.$$
As a first example of lattice extensions, we demonstrate the following construction of a small-determinant extension of a sublattice inside of the integer lattice~$\zed^n$. 

\begin{thm} \label{det_ext} Let $\bx_1,\dots,\bx_k$ be linearly independent vectors in $\zed^n$ and let
$$\Omega = \spn_{\zed} \left\{ \bx_1,\dots,\bx_k \right\} \subset \zed^n$$
be the sublattice of rank $k$ spanned by these vectors. Then there exists a full-rank extension $\Omega' \subseteq \zed^n$ of $\Omega$ so that
$$\det \Omega' = \gcd( \bx_1 \wedge \dots \wedge \bx_k ).$$
Further, if $k=n-1$ then there exists $\bwy \in \zed^n$ so that $\Omega' = \spn_{\zed} \left\{ \Omega, \bwy \right\}$ and
$$\|\bwy\| \leq \left\{ \left( \frac{\gcd( \bx_1 \wedge \dots \wedge \bx_k )}{\det \Omega} \right)^2 + \mu(\Omega)^2 \right\}^{1/2}.$$
\end{thm}

\noindent
Throughout this paper, the wedge product of vectors, $\bx_1 \wedge \dots \wedge \bx_k$ as above in the Grassmann algebra, is identified with the corresponding vector of Grassmann-Pl\"ucker coordinates; see Chapter 1 of~\cite{schmidt1991} for details. We prove Theorem~\ref{det_ext} in Section~\ref{Z-ext}, where we also explain how it can be generalized to any lattice in~$\real^n$ (Remark~\ref{det_ext_gen}). Lattice extensions play an implicit important role in a variety of contexts, for instance in lattice packing and covering constructions such as lamination (see \cite{conway:sloane}, \cite{martinet}), in reduction theory when constructing Minkowski or HKZ reduced bases (see \cite{lek}, \cite{martinet}), as well as in constructions of primitive collections in a lattice (see \cite{max_lenny}). Further, the idea of lattice extensions was recently used to construct a family of counterexamples to the famous covering conjecture of Woods (see~\cite{regev1}). This being said, we have not seen lattice extensions studied explicitly. Our main goal in this note is to explore lattice extensions with control over their geometric invariants. In particular, we say that $\Lambda$ is a {\it successive minima extension} of $L$ if $\Lambda$ is an extension of $L$ such that
$$\lambda_j(\Lambda) = \lambda_j(L)\ \forall\ 1 \leq j \leq k.$$
Also, we call $\Lambda$ an {\it equal covering extension} of $L$ if $\Lambda$ is an extension of $L$ such that
$$\mu(\Lambda) = \mu(L).$$
Given a lattice $L$ of rank $k < n$ in $\real^n$, it is easy to construct a rank-$(k+1)$ successive minima extension $\Lambda$ of $L$: we can simply take $\bu \in \real^n$ to be a vector perpendicular to $\spn_{\real} L$ of norm $ > \lambda_k(L)$ and define $\Lambda = \spn_{\zed} \{ L, \bu \}$. It is a more delicate problem to construct such an extension inside of a given full-rank lattice in~$\real^n$, since such a perpendicular vector $\bu$ may simply not exist inside of our given lattice. Our next result addresses this problem while controlling the $(k+1)$-st successive minimum of the constructed extension.

\begin{thm} \label{sm_ext} Let $\Lambda \subset \real^n$ be a lattice of full rank, and let $L_k \subset \Lambda$ be a sublattice of rank $1 \leq k < n$. There exists a sublattice $L_{k+1} \subset \Lambda$ of rank $k+1$ such that $L_k \subset L_{k+1}$ is a lattice extension, $\lambda_j(L_{k+1}) = \lambda_j(L_k)$ for all $1 \leq j \leq k$ and
\begin{equation}
\label{l_bnd}
\lambda_{k+1}(L_{k+1}) \leq \frac{\lambda_k(L_k)(v_*^2 + \sqrt{1-v_*^2})}{\sqrt{1-v_*^4}} + 2\mu,
\end{equation}
where $\mu$ is the covering radius of $\Lambda$ and $v_*$ is the smallest root of the polynomial
$$p(v) = \left( \frac{\mu^2}{\lambda_k^2} (1 - v^4) - v^2 (v^4 - v^2 + 1) \right)^2 - \left( \frac{2 \mu^2}{\lambda_k^2} v (1-v^4) + 2v^4 \right)^2 (1 - v^2)$$
in the interval~$(0,1)$: such $v_*$ necessarily exists.
\end{thm}

\noindent
We prove Theorem~\ref{sm_ext} in Section~\ref{succ_min}. We also include an alternate version of the bound for Theorem~\ref{sm_ext} suggested to us by one of the referees (Remark~\ref{ref_rem}). The situation is more complicated with equal covering extensions: they do not necessarily exist inside of a given lattice. Our next result is a full characterization of planar lattices that are equal covering extensions of some lattice of rank one. Let $\be_1 = \begin{pmatrix} 1 \\ 0 \end{pmatrix} \in \real^2$ and $E_1 = \zed \be_1 \subset \real^2$ be a lattice of rank one in the plane. Then the covering radius of $E_1$ is $\mu(E_1) = 1/2$. More generally, for a rank-one lattice $\zed \bu \in \real^n$ its covering radius is $\frac{1}{2} \|\bu\|$.

\begin{thm} \label{cover_ext} A lattice $\Lambda \subset \real^2$ is an equal covering extension of $E_1$ if and only if
\begin{equation}
\label{cv_ext}
\Lambda = \Lambda(\alpha) := \begin{pmatrix} \alpha & \alpha - 1 \\ \sqrt{\alpha-\alpha^2} & \sqrt{\alpha-\alpha^2} \end{pmatrix} \zed^2
\end{equation}
for some real number $0 < \alpha < 1$. More generally, a lattice $\Lambda \subset \real^n$ of rank~$2$ is an equal covering extension of a rank-one lattice $L \subset \Lambda$ if and only if it is isometric to some lattice of the form $\det(L) \Lambda(\alpha)$, where $\Lambda(\alpha)$ is as in~\eqref{cv_ext}. 
\end{thm}

\noindent
We discuss covering radii of planar lattices and prove Theorem~\ref{cover_ext} with some corollaries in Section~\ref{cover}. In particular, we show that the lattice coming from the ring of integers of a quadratic number field $\que(\sqrt{D})$, for a squarefree rational integer $D$, via Minkowski embedding into~$\real^2$ is an equal covering extension of a rank-one sublattice if and only if $D \not\equiv 1\ (\md 4)$. We also construct orthogonal equal covering extensions in any dimension, proving the following result.

\begin{thm} \label{mu_orth} Let \(\Lambda_k \subset \real^n\) be an orthogonal lattice of rank \(k < n\). There exists an orthogonal lattice \(\Lambda_{k+1} \subset \real^n\) of rank \(k+1\) so that \(\Lambda_k \subset \Lambda_{k+1}\)  is a lattice extension and \(\mu(\Lambda_{k+1}) = \mu(\Lambda_k)\). Further, if \(\bz\) is a deep hole of \(\Lambda_k\) it is also a deep hole of \(\Lambda_{k+1}\).
\end{thm}

\noindent
Recall that, given a lattice $\Lambda \subset \real^n$ a vector $\bz \in \spn_{\real} \Lambda$ is called a {\it deep hole} of $\Lambda$ if it is furthest from the lattice, i.e.
$$d(\bz,\Lambda) = \max \left\{ d(\bwy, \Lambda) : \bwy \in \spn_{\real} \Lambda \right\},$$
where $d(\bwy, \Lambda) := \min \{ \|\bx - \bwy\| : \bx \in \Lambda \}$. Thus the covering radius of the lattice is the distance from the origin to the nearest deep hole. We discuss deep holes of lattices in some detail in Section~\ref{deep_holes} with a special focus on the two-dimensional situation. In particular, we obtain necessary and sufficient conditions for the deep holes of a lattice $\Lambda \subset \real^2$ to have finite order as elements of the torus quotient group $\real^2/\Lambda$ and give a bound on this order (Theorem~\ref{dp_hole_fin}).

Before we proceed, let us recall a few more standard notions of lattice theory. The isometries of a Euclidean space are given by real orthogonal matrices, and two lattices are {\it isometric} if there exists an isometry taking one to the other. Two lattices are {\it similar} if their scalar multiples are isometric for some choice of scalars. Both, isometry and similarity are equivalence relations on lattices of the same rank. A lattice is called {\it well-rounded} (abbreviated WR) if all of its successive minima are equal; this property is preserved under similarity. 

\bigskip

\section{Small-determinant extensions in $\zed^n$}
\label{Z-ext}

In this section we present the first example of lattice extensions, proving Theorem~\ref{det_ext}. Let $\bx_1,\dots,\bx_k$ be $k$ vectors in $\real^n$, $1 \leq k < n$. As we mentioned above, the wedge product $\bx_1 \wedge \dots \wedge \bx_k$ can be identified with the vector of Pl\"ucker coordinates in~$\real^{\binom{n}{k}}$, i.e. determinants of $k \times k$ submatrices of the $n \times k$ matrix $(\bx_1\ \dots\ \bx_k)$. 

\begin{prop} \label{Zext-1} Let $\bx_1,\dots,\bx_k$ be linearly independent vectors in $\zed^n$ and let
$$\Omega = \spn_{\zed} \left\{ \bx_1,\dots,\bx_k \right\} \subset \zed^n$$
be the sublattice of rank $m$ spanned by these vectors. Then there exists a full-rank extension $\Omega' \subseteq \zed^n$ of $\Omega$ so that
$$\det \Omega' = \gcd( \bx_1 \wedge \dots \wedge \bx_k ).$$
\end{prop}

\proof
Let $\bar{\Omega} = \zed^n \cap \spn_{\real} \Omega$, then $\bar{\Omega} \subset \zed^n$ is a sublattice of rank $k$ containing $\Omega$ such that $\zed^n/\bar{\Omega}$ is torsion free. Hence any basis of $\bar{\Omega}$ is extendable to $\zed^n$. Let $\bwy_1,\dots,\bwy_k$ be a basis for $\bar{\Omega}$ extended to a basis for $\zed^n$ by $\bwy_{k+1},\dots,\bwy_n$. Since $\bx_1,\dots,\bx_k$ and $\bwy_1,\dots,\bwy_k$ are two collections of integer vectors spanning the same subspace of $\real^n$, the vectors of Pl\"ucker coordinates represent the same rational projective point. Further, since  the collection $\bwy_1,\dots,\bwy_k$  is extendable to a basis of $\zed^n$, the Pl\"ucker coordinates of this collection must be relatively prime (Lemma~2 on p.15 of~\cite{cassels}). Hence
$$\bx_1 \wedge \dots \wedge \bx_k = g \left( \bwy_1 \wedge \dots \wedge \bwy_k \right)$$
for some integer $g$, and so $g = \gcd( \bx_1 \wedge \dots \wedge \bx_k )$. Define 
$$\Omega' = \spn_{\zed} \left\{ \bx_1,\dots,\bx_k,\bwy_{k+1},\dots,\bwy_n \right\}.$$
By the bi-linearity of the wedge product, 
$$\det \Omega' = \bx_1 \wedge \dots \wedge \bx_k \wedge \bwy_{k+1} \wedge \dots \wedge \bwy_n = g \left( \bwy_1 \wedge \dots \wedge \bwy_k \wedge \bwy_{k+1} \wedge \dots \wedge \bwy_n \right),$$
and since $\bwy_1 \wedge \dots \wedge \bwy_k \wedge \bwy_{k+1} \wedge \dots \wedge \bwy_n = \det \zed^n = 1$, we have that $\det \Omega' = g$.
\endproof

\begin{cor} \label{Zext-2} Let the notation be as in Proposition~\ref{Zext-1} with $k=n-1$. Then there exists $\bwy \in \zed^n$ so that $\Omega' = \spn_{\zed} \left\{ \Omega, \bwy \right\}$ and
$$\|\bwy\| \leq \left\{ \left( \frac{\gcd( \bx_1 \wedge \dots \wedge \bx_k )}{\det \Omega} \right)^2 + \mu(\Omega)^2 \right\}^{1/2}.$$
\end{cor}

\proof
Write $A = \begin{pmatrix} \bx_1 & \dots & \bx_{n-1} \end{pmatrix}$ for the corresponding basis matrix of $\Omega$ and let $\Omega'$ be as given by Proposition~\ref{Zext-1}. This means that there exists $\bz \in \zed^n$ such that $\Omega' = \spn_{\zed} \{ \Omega, \bz \}$, so $\det \Omega' = \gcd( \bx_1 \wedge \dots \wedge \bx_m )$. Let $\rho_{\Omega} = A(A^{\top} A)^{-1}A^{\top}$ be the orthogonal projection onto $\spn_{\real} \Omega$. Let
$$\PP = \left\{ \sum_{i=1}^{n-1} a_i \bx_i : 0 \leq a_i < 1\ \forall\ 1 \leq i \leq n-1 \right\},\ \PP' = \left\{ \bu + a\bz : \bu \in \PP,\ 0 \leq a < 1 \right\}$$
be fundamental parallelepipeds for $\Omega$ and $\Omega'$, respectively. Then
\begin{eqnarray*}
\gcd( \bx_1 \wedge \dots \wedge \bx_k ) & = & \det \Omega' = \Vol_n(\PP') \\
& = & \Vol_{n-1}(\PP) \left\| (I_n - \rho_{\Omega}) \bz \right\| = \det \Omega \left\| (I_n - \rho_{\Omega}) \bz \right\|,
\end{eqnarray*}
hence
$$\left\| (I_n - \rho_{\Omega}) \bz \right\| = \frac{\gcd( \bx_1 \wedge \dots \wedge \bx_k )}{\det \Omega}.$$
On the other hand, $\rho_{\Omega} \bz \in \spn_{\real} \Omega$, and by definition of the covering radius $\mu$ of $\Omega$, there exists $\bv \in \Omega$ such that $\|\rho_{\Omega} \bz - \bv \| \leq \mu$. Let $\bwy = \bz - \bv$, then $\bwy \in \zed^n$ and 
$$\rho_{\Omega} \bwy = \rho_{\Omega} \bz - \rho_{\Omega} \bv = \rho_{\Omega} \bz - \bv,$$
since $\bv \in \spn_{\real} \Omega$. Then $(I_n - \rho_{\Omega}) \bwy = (I_n - \rho_{\Omega}) \bz$ and
$$\Omega' = \spn_{\zed} \left\{ \Omega, \bz \right\} = \spn_{\zed} \left\{ \Omega, \bwy \right\}.$$
Therefore, by Pythagorean theorem,
\begin{eqnarray*}
\|\bwy\|^2 & = & \left\| (I_n - \rho_{\Omega}) \bwy \right\|^2 + \| \rho_{\Omega} \bwy \|^2 = \left\| (I_n - \rho_{\Omega}) \bz \right\|^2 + \| \rho_{\Omega} \bz - \bv \|^2 \\
& \leq & \left( \frac{\gcd( \bx_1 \wedge \dots \wedge \bx_k )}{\det \Omega} \right)^2 + \mu^2.
\end{eqnarray*}
The result follows.
\endproof

Now Theorem~\ref{det_ext} follows by combining Proposition~\ref{Zext-1} with Corollary~\ref{Zext-2}.

\begin{rem} \label{det_ext_gen} Let $\Lambda = A\zed^m \subset \real^n$ be a lattice of rank $m \leq n$ and let $\bz_1,\dots,\bz_k$, $k \leq m$, be linearly independent vectors in $\Lambda$. Then for each $1 \leq i \leq k$, $\bz_i = A \bx_i$, where $\bx_1,\dots,\bx_k \in \zed^n$ are also linearly independent. Let
$$\Omega = \spn_{\zed} \left\{ \bx_1,\dots,\bx_k \right\} \subset \zed^m$$
be the sublattice of rank $k$ spanned by these vectors and let $\Omega'$ be an extension of $\Omega$ in $\zed^n$ guaranteed by Proposition~\ref{Zext-1}. Then $A\Omega = \spn_{\zed} \left\{ \bz_1,\dots,\bz_k \right\} \subseteq \Lambda$ and $A\Omega' \subseteq \Lambda$ is an extension of $A\Omega$ with
$$\det A\Omega' = \sqrt{ \det (A^{\top} A) }\ \det \Omega' = \det \Lambda \det \Omega'.$$
Further, if $k=m-1$ then there exists $\bwy \in \zed^m$ so that $A\Omega' = \spn_{\zed} \left\{ A\Omega, A\bwy \right\}$ and $\|\bwy\|$ bounded as in Corollary~\ref{Zext-2}.
\end{rem}

\bigskip

\section{Successive minima extensions}
\label{succ_min}

In this section, we prove Theorem~\ref{sm_ext}. We want to construct a sublattice $L_{k+1} \subset \Lambda$ of rank $k+1$ such that $L_k \subset L_{k+1}$, $\lambda_j(L_{k+1}) = \lambda_j(L_k)$ for all $1 \leq j \leq k$ and $\lambda_{k+1}(L_{k+1})$ is as small as possible. To prove the theorem, we first need an auxiliary lemma. Write $\lambda_1,\dots,\lambda_k$ for the successive minima of $L_k$ and let $V_k = \spn_{\real} L_k$, $\theta \in (0,\pi/2]$, and
\begin{equation}
\label{cone}
C_{\theta}(V_k) = \left\{ \bx \in \real^n : \a(\bx,\bwy) \in [\theta,\pi-\theta]\ \forall\ \bwy \in V_k \right\},
\end{equation}
where $\a(\bx,\bwy)$ stands for the angle between two vectors.

\begin{lem} \label{nrm_min} If $\bx \in C_{\theta}(V_k)$ and
$$\|\bx\| \geq \frac{\lambda_k (\cot \theta \cos \theta + 1)}{\sqrt{1+\cos^2 \theta}},$$
then $\|\bx + \bwy\| \geq \lambda_k$ for every $\bwy \in V_k$.
\end{lem}

\proof
For $\bx \in C_{\theta}(V_k)$ and $\bwy \in V_k$, define
$$f(\bx,\bwy) = \| \bx + \bwy\|^2 = \|\bx\|^2 + \|\bwy\|^2 + 2 \|\bx\| \|\bwy\| \cos \a(\bx,\bwy).$$
We want to guarantee that $f(\bx,\bwy) \geq \lambda_k^2$ for all $\bwy \in V_k$. Let us write $t = \|\bx\|$, $z = \|\bwy\|$, and notice that
$$f(\bx,\bwy) - \lambda_k^2 \geq g(t,z) := t^2 + z^2 - 2 tz \cos \theta - \lambda_k^2,$$
thus we want to find a lower bound on $t$ that would guarantee $g(t,z) \geq 0$ for all $z > 0$. In other words, we want
$$t \geq h(z) := z \cos \theta + \sqrt{ \lambda_k^2 - z^2 \sin^2 \theta}$$
for all $z > 0$. Notice that $h(z)$ is real-valued if and only if $z \leq \frac{\lambda_k}{\sin \theta}$, then let us find the value of $z$ that maximizes $h(z)$. Differentiating $h(z)$ and setting the derivative equal to zero, we obtain
$$z_* = \frac{\lambda_k \cot \theta}{\sqrt{1+\cos^2 \theta}},$$
the point at which $h(z)$ assumes its maximum value of
$$h(z_*) = \frac{\lambda_k (\cot \theta \cos \theta + 1)}{\sqrt{1+\cos^2 \theta}}.$$
Thus, taking $t = \|\bx\|$ to be $\geq$ than this value ensures that $\|\bx + \bwy\| \geq \lambda_k$ for every $\bwy \in V_k$, as required.
\endproof

\proof[Proof of Theorem~\ref{sm_ext}]
Let us write $B_n(r)$ for the ball of radius $r > 0$ centered at the origin in~$\real^n$. Let $\theta \in (0,\pi/2]$ and
\begin{equation}
\label{rth}
r(\theta) = \frac{\lambda_k (\cot \theta \cos \theta + 1)}{\sqrt{1+\cos^2 \theta}}.
\end{equation}
Then Lemma~\ref{nrm_min} guarantees that for any vector $\bx \in \Lambda \cap \left( C_{\theta}(V_k) \setminus B_n(r(\theta)) \right)$ the lattice $M = \spn_{\zed} \left\{ L_k, \bx \right\}$ satisfies $\lambda_j(M) = \lambda_j(L_k)$ for all $1 \leq j \leq k$ and $\lambda_{k+1}(M) \leq \|\bx\|$. Hence we want to minimize
$$\lambda_{k+1}(\theta) := \min \left\{ \|\bx\| : \bx \in \Lambda \cap (C_{\theta}(V_k) \setminus B_n(r(\theta)) \right\}$$
as a function of $\theta$. 

\begin{figure}
\centering
\includegraphics[scale=0.4]{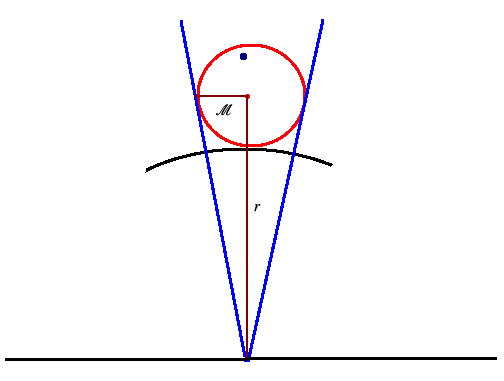}
\caption{Cone construction with the lattice point (in blue) caught in the ball (red) of covering radius.}\label{cone1}
\end{figure}

Let $\mu$ be the covering radius of $\Lambda$, then any translated copy of the ball of radius $\mu$ in $\real^n$ must contain a point of $\Lambda$. Suppose that $\theta \in (0,\pi/2]$ is such that
\begin{equation}
\label{ball_theta}
B'_n(\mu) \subset \left( C_{\theta}(V_k) \cap B_n(r(\theta)+2\mu) \right) \setminus B_n(r(\theta)),
\end{equation}
where $B'_n(\mu)$ is such a translated copy. Then $C_{\theta}(V_k) \setminus B_n(r(\theta))$ would be guaranteed to contain a point $\bx$ of $\Lambda$ with
\begin{equation}
\label{x_bnd}
\|\bx\| \leq r(\theta) + 2\mu,
\end{equation}
so that we can take $L_{k+1} = \spn_{\zed} \left\{ L_k,\bx \right\}$. Notice that $\theta$ satisfying condition~\eqref{ball_theta} exists. Indeed, for any $\theta$ the set $B_n(r(\theta)+2\mu) \setminus B_n(r(\theta))$ contains a ball of radius $\mu$, and hence $\theta$ can always be chosen small enough so that the cone $C_{\theta}(V_k)$ is sufficiently wide to contain this ball. In fact, we can choose $\theta$ so that the line segment from $\bo$ to the center of this ball $B'_n(\mu)$ has length $r(\theta)+\mu$ and makes the angle $\pi/2-\theta$ with any line in the boundary of $C_{\theta}(V_k)$ emanating from the center and tangent to the ball $B'_n(\mu)$. These conditions result in a right triangle with legs $r(\theta)+\mu$ and $\mu$ and the angle $\pi/2-\theta$ opposite to the second leg (see Figure~\ref{cone1} for a graphical illustration of this argument). Hence we have the equation
$$\tan(\pi/2-\theta) = \frac{\mu}{r(\theta)+\mu}.$$

\noindent
Using \eqref{rth}, along with the fact that $\tan(\pi/2-\theta) = \cot \theta$, writing $v = \cos \theta$ and simplifying, we obtain the following relation in terms of $v$:
$$\mu \left( \sqrt{1-v^2} - v \right) = \frac{\lambda_k \left( v^2 + \sqrt{1-v^2} \right) v}{\sqrt{1-v^4}},$$
which transforms into the following polynomial equation:
\begin{equation}	
\label{v_eq}
\left( \frac{\mu^2}{\lambda_k^2} (1 - v^4) - v^2 (v^4 - v^2 + 1) \right)^2 = \left( \frac{2 \mu^2}{\lambda_k^2} v (1-v^4) + 2v^4 \right)^2 (1 - v^2).
\end{equation}
It follows from our construction that this equation has at least one solution $v$ in the interval $(0,1)$. Then $r(\theta)$ as a function of $v$ becomes
$$r(v) = \frac{\lambda_k(v^2 + \sqrt{1-v^2})}{\sqrt{1-v^4}},$$
which is an increasing function of $v$ in the interval $(0,1)$. Hence, to minimize the bound~\eqref{x_bnd}, we can pick $v_*$ to be the smallest root of the equation~\eqref{v_eq} in the interval~$(0,1)$. In other words, we are maximizing our choice of the angle $\theta$ for which the condition~\eqref{ball_theta} holds. The inequality~\eqref{l_bnd} follows.
\endproof
\smallskip

\begin{rem} \label{ref_rem} We also present here an alternate bound for the $(k+1)$-st successive minimum of an extension lattice $L_{k+1}$ of $L_k$ in $\Lambda$ with $\lambda_j(L_{k+1}) = \lambda_j(L_k)$ for all $1 \leq j \leq k$ as in Theorem~\ref{sm_ext}, as suggested to us by one of the anonymous referees:
\begin{equation}
\label{rr1}
\lambda_{k+1}(L_{k+1}) \leq \left( \max \{ \lambda_{k+1}(\Lambda), 2\lambda_k(L_k) \}^2 + \mu(L_k)^2 \right)^{1/2}.
\end{equation}
To prove this bound, let $\bu_1,\dots,\bu_{k+1}$ be linearly independent vectors in $\Lambda$ corresponding to the successive minima $\lambda_1(\Lambda),\dots,\lambda_{k+1}(\Lambda)$, respectively. Since $\dim V_k = k$, there must exist some $\bu_i$ among these vectors which is not in $V_k$. Let $\bx$ be the projection of $\bu_i$ into $V_k^{\perp}$, the orthogonal complement of $V_k$. First assume that $\|\bx\| \geq \lambda_k(L_k)$, then $\|\bx+\bwy\| \geq \lambda_k(L_k)$ for any $\bwy \in V_k$. To see this, notice that
\begin{equation}
\label{V_k-orth}
\|\bx+\bwy\|^2 = \|\bx\|^2 + 2 \bx^{\top} \bwy + \|\bwy\|^2 = \|\bx\|^2 + \|\bwy\|^2 \geq \|\bx\|^2,
\end{equation}
since $\bx \in V_K^{\perp}$, $\bwy \in V_k$, and so $\bx^{\top} \bwy = 0$.
The translated subspace $\bx + V_k$ contains an affine copy of the lattice $L_k$, and hence there must exist a point of the lattice~$\Lambda$, call it $\bx_{k+1}$, in the set $\bx+B_{V_k} \left( \mu(L_k) \right)$, where $B_{V_k} \left( \mu(L_k) \right)$ is the ball of radius $\mu(L_k)$ centered at the origin in $V_k$. Set $L_{k+1} = \spn_{\zed} \{ L_k, \bx_{k+1} \}$, then~\eqref{V_k-orth} implies that for every $\bz \in L_{k+1} \setminus L_k$, $\|\bz\| \geq \|\bx\| \geq \lambda_k(L_k)$. Hence, $\lambda_j(L_{k+1}) = \lambda_j(L_k)$ for all $1 \leq j \leq k$ and
$$\lambda_{k+1}(L_{k+1})^2 \leq \|\bx_{k+1}\|^2 \leq \|\bx\|^2 + \mu(L_k)^2 \leq \|\bu_i\|^2 + \mu(L_k)^2 \leq \lambda_{k+1}(\Lambda)^2 + \mu(L_k)^2,$$
implying~\eqref{rr1}. On the other hand, suppose that $\|\bx\| < \lambda_k(L_k)$. We can then replace $\bx$ by its integer multiple $\bx'$, chosen in such a way that $\lambda_k(L_k) \leq \|\bx'\| \leq 2\lambda_k(L_k)$, and proceed as above. We thank the referee for suggesting this elegant argument.
\end{rem}

\begin{rem} \label{ss_comp} We also want to discuss Theorem~\ref{sm_ext} in the general context of reduction theory. This result can be loosely compared to the construction of a {\it canonical filtration} of a lattice as originally defined by Grayson and Stuhler (see Casselman's survey paper~\cite{casselman} for a detailed discussion, as well as~\cite{regev2} for a recent application of the canonical filtration). This is a unique flag of sublattices
$$\{\bo\} = \Lambda_0 \subset \Lambda_1 \subset \dots \subset \Lambda_n = \Lambda$$
in a lattice $\Lambda$ such that $\rk (\Lambda_k) = k$ and $\det (\Lambda_n)^{1/n} > \det (\Lambda_k)^{1/k}$ for every $k < n$, where
$$\det (\Lambda_k) = \min \left\{ \det (\Omega) : \Omega \subset \Lambda,\ \rk (\Omega) = k \right\}.$$
A lattice $\Lambda$ is called {\it semi-stable} if the canonical filtration is $\Lambda_0 \subset \Lambda_n$, i.e. if for each sublattice $\Omega \subseteq \Lambda$,
\begin{equation}
\label{ss_det}
\det (\Lambda)^{1/\rk(\Lambda)} \leq \det (\Omega)^{1/\rk(\Omega)}.
\end{equation}
This family of lattices is important in reduction theory. Y. Andre explains in~\cite{andre}:
\smallskip

{\it Reduction theory aims at estimating the length of short vectors, and more generally the (co)volumes of small sublattices of lower ranks, of lattices of given rank and (co)volume, and at combining lower and upper bounds to get finiteness results. A better grasp on lower bounds comes from the more recent part of reduction theory which deals with semistability and slope filtrations (heuristically, semistability means that the Minkowski successive minima are not far from each other, cf.~\cite{borek})}
\smallskip

\noindent
On the other hand, our Theorems~\ref{det_ext} and~\ref{sm_ext} give constructions of lattice extensions of a given sublattice within an ambient lattice with small determinant (= (co)volume) and successive minima, respectively, while preserving the geometric properties of the sublattice that is being extended.

Additionally, if we start with a rank-one sublattice $L_1 \subset \Lambda$ and recursively apply the construction described in the proof of Theorem~\ref{sm_ext}] to obtain a sublattice $L_n \subseteq \Lambda$ of full rank, the collection of vectors we build to bound the successive minima at every step will be a basis, call this basis $\bx_1,\dots,\bx_n$. We can choose the angle $\theta_k$ between $\bx_{k+1}$ and $V_k = \spn_{\real} L_k = \spn_{\real} \{ \bx_1,\dots,\bx_k \}$ for each $1 \leq k \leq n-1$ to be in the interval~$[\pi/3,\pi/2]$ instead of the interval $(0,\pi/2]$ we used in the proof of Theorem~\ref{sm_ext}]: this is always possible at the expense of the $(k+1)$-st successive minimum $\lambda_{k+1}$ being larger, since the cone $C_{\theta_k}(V_k)$ eventually becomes wide enough to contain a ball of radius $\mu$. Then we can ensure that the resulting lattice $L_n$ is weakly nearly orthogonal in the sense of~\cite{baraniuk} and~\cite{lf_dk}. Specifically, a lattice is called weakly nearly orthogonal if it contains an ordered basis with the angles between the $(k+1)$-st basis vector and the subspace spanned by the previous $k$ falling in the interval $[\pi/3,\pi/2]$. Weakly nearly orthogonal lattices have applications in image compression and digital communications. 
\end{rem}

\bigskip

\section{Deep holes of planar lattices}
\label{deep_holes}

We start this section with the following simple but useful technical lemma.

\begin{lem} \label{simplex} Let $\bx_1,\ldots, \bx_m$ be linearly independent points in $\real^n$, $m \leq n$. There exist points \(\bz \in \real^n\) so that \(\|\bz\| = \|\bz-\bx_i\|\) for all \(1 \leq i\leq m\), and these points are solutions to 
\begin{equation}\label{simplex-eq}
  \begin{pmatrix}
  \bx_1^\top\\\vdots\\\bx_m^\top
  \end{pmatrix}
\bz
  =\frac{1}{2} 
  \begin{pmatrix}
  \|\bx_1\|^2 \\ \vdots \\ \|\bx_m\|^2
  \end{pmatrix}.
\end{equation} 
\end{lem}

\proof
If \(\bz\) is equidistant from \(\bx_i\) and \(\textbf{0}\) then \(\bz\) lies in the hyperplane orthogonal to \(\bx_i\) that passes through the point \(\bx_i/2\). That is 
\[
\bz \cdot \bx_i = \text{proj}_{\bx_i}(\bz) \cdot \bx_i = \frac{\|\bx_i\|^2}{2}.
\]
Since this is true for each \(1 \leq i \leq n\) this gives the linear system in \eqref{simplex-eq}.
\begin{figure}
    \centering
    \begin{tikzpicture}
        \draw[thick, black] (-1,0) -- (4,0);
        \draw[thick, black] (0,-1) -- (0,4);
        \draw[->] (0,0)--(3,1);
        \filldraw (3,1) node[anchor = west]{$\bx_1$};
        \draw[->] (0,0)--(1,3);
        \filldraw (1,3) node[anchor = west]{$\bx_2$};
        \draw[dotted] (-1,2) -- (4,0.33);
        \draw[dotted] (2,-1) -- (0.33,4);
        \filldraw (1.25, 1.25) circle (2pt) node[anchor=west]{$\bz$};
        \filldraw (1.5, 0.5) circle (2pt) node[anchor=north west]{$\text{proj}_{\bx_1}(\bz)$};
    \end{tikzpicture}
\caption{Point $\bz$, equidistant from $\bo,\bx_1,\bx_2$, to illustrate the proof of Lemma~\ref{simplex}.}\label{proj_fig}
\end{figure}
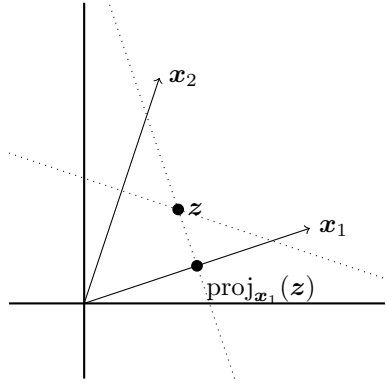
\endproof

Our main goal here is to describe some properties of the deep holes of lattices, focusing especially on the two-dimensional case. Our first basic observation is that if $\bz$ is a deep hole of a lattice $\Lambda \subset \real^n$, then so is $-\bz$: this follows by the fact that $-\Lambda = \Lambda$, since lattices are symmetric about the origin. Recall that every $2$-dimensional lattice has a basis consisting of vectors corresponding to successive minima, and such a basis can always be chosen so that the angle between the basis vectors is in the interval $[\pi/3,\pi/2]$; we call this basis a {\it minimal basis} for the lattice (see, e.g., \cite{hex} for details). Then we have the following observation.

\begin{lem} \label{deep_hole} Let $\Lambda \subset \real^2$ be a lattice of rank~$2$ with minimal basis $\bx,\bwy$ and angle $\theta \in [\pi/3,\pi/2]$ between these basis vectors. Write $\lambda_1,\lambda_2$ for the successive minima of $\Lambda$, so that $0 < \lambda_1 = \|\bx\| \leq \lambda_2 = \|\bwy\|$. Then the fundamental parallelogram
$$\PP = \left\{ s \bx + t \bwy : 0 \leq s,t < 1 \right\}$$
contains two deep holes $\bz_1,\bz_2$ with $\bz_1+\bz_2 \in \Lambda$. If the angle $\theta = \pi/2$, then $\bz_1 = \bz_2$ is the center of $\PP$, and we say that this deep hole has multiplicity~$2$.
\end{lem}

\begin{figure}
\centering
\includegraphics[scale=0.4]{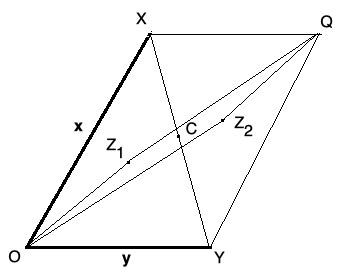}
\caption{Fundamental parallelogram $\PP$ of $\Lambda$ with deep holes $\bz_1$ and $\bz_2$.}\label{pgrm}
\end{figure}

\proof
Let us label the vertices of $\PP$ as follows: $O$ for the origin, $X$ for the endpoint of the vector $\bx$, $Y$ for the endpoint of the vector $\bwy$, and $Q$ for the endpoint of the vector $\bx+\bwy$. The parallelogram $\PP$ can be split into two congruent triangles: $OXY$ and $QYX$. Then the endpoints of the deep holes of $\Lambda$ contained in $\PP$ are the centers of the circles circumscribed around these two triangles, call them $Z_1$ and $Z_2$, respectively, and let $\bz_1,\bz_2$ be vectors with the endpoints $Z_1,Z_2$ (see~\cite{cover-2}). The two triangles are symmetric to each other about the center $C$ of $\PP$, which means that reflection with respect to $C$ maps the line segment $OZ_1$ onto the line segment $QZ_2$. This means that $OZ_1QZ_2$ is a parallelogram with $OQ$ as its longer diagonal, and hence the corresponding vector is the sum $\bz_1+\bz_2$. Since its endpoint is $Q$, a vertex of $\PP$, this vector is in~$\Lambda$ (see Figure~\ref{pgrm} for a graphical illustration of this argument). If $\theta = \pi/2$, then the deep hole of each of the triangles is in the center of the hypothenuse of its corresponding right triangle, i.e. at the center point $C$ of $\PP$; in this case, the two deep holes coincide, so $Z_1 = Z_2 = C$. 
\endproof

An immediate implication of Lemma~\ref{deep_hole} is that deep holes $\bz_1,\bz_2$ are each other's inverses in the additive abelian group $\real^2/\Lambda$. Further, $\bz_1$ is an element of order two in this group if and only if the angle $\theta = \pi/2$; in this case $\bz_1=\bz_2$. On the other hand, $\bz_1,\bz_2$ can be elements of finite order in other situations too. For instance, in the hexagonal lattice
$$L_{\pi/3} = \begin{pmatrix} 1 & \frac{1}{2} \\ 0 & \frac{\sqrt{3}}{2} \end{pmatrix} \zed^2$$
the deep holes are $\bz_1 = (1/2, 1/(2\sqrt{3}))$, $\bz_2 = (1, 1/\sqrt{3})$ have order three in the group $\real^2/L_{\pi/3}$, while the lattice 
$$L' = \begin{pmatrix} 1 & \frac{1}{2} \\ 0 & \sqrt{3} \end{pmatrix} \zed^2$$
has a deep hole $\bz_1 = (1/2, 11\sqrt{3}/24)$ satisfying the condition
$$48 \bz_1 = 13 (1,0) + 22 (1/2, \sqrt{3}) \in L',$$
which makes $\bz_1$ an element of order dividing 48 in the group $\real^2/L'$. These observations raise a natural question: when does a deep hole of $\Lambda \subset \real^2$ have finite order as an element of the group $\real^2/\Lambda$? 

\begin{thm} \label{dp_hole_fin} Let $\Lambda \subset \real^2$ be a full-rank lattice with successive minima $\lambda_1,\lambda_2$ and corresponding minimal basis vectors $\bx_1, \bx_2$. A deep hole $\bz$ of $\Lambda$ has finite order in the group \(\real^2/\Lambda\) if and only if \(\Lambda\) is orthogonal or there exist rational numbers \(p,q \) so that \(p\lambda_1^2 = \bx_1\cdot \bx_2 = q\lambda_2^2\). Moreover, if \(\lambda_1^2, \lambda_2^2,\bx_1\cdot \bx_2\in \zed\) then the order of \(\bz\) in \(\real^2/\Lambda\) is \(\leq 12 \sqrt{3}\ \lambda_2^4\).
\end{thm}

\proof
As we discussed above, if \(\Lambda\) is orthogonal then the deep hole always has order \(2\) in \(\real^2/\Lambda\), hence we assume \(\Lambda\) is not orthogonal. As we indicated above, we can assume that the minimal basis vectors \(\bx_1, \bx_2\) are chosen so that the angle $\theta$ between them satisfies \(\pi/3 \leq \theta \leq \pi/2\). If $\bz$ is the equidistant from \(\bx_1,\bx_2\) and the origin then \({\bz}\) is a deep hole of \(\Lambda\) and is contained in the convex hull of \(\{0,\bx_1,\bx_2\}\).
By Lemma~\ref{simplex},
\begin{equation}\label{order-rat1}
\bz \cdot \bx_1 = \frac{\lambda_1^2}{2},\ \bz \cdot \bx_2 =  \frac{\lambda_2^2}{2}.
\end{equation}
Now suppose that $\bz$ has finite order in \(\real^2/\Lambda\). Then there are integers $a,b,c$ so that $c \neq 0$ and
$$a\bx_1 + b\bx_2 = c \bz.$$
In fact, the pairs \(\bz, \bx_1\) and \(\bz, \bx_2\) are linearly independent so \(a,b,c\) are all nonzero. Taking scalar products of both sides of this equation with $\bx_1$ and $\bx_2$, and applying~\eqref{order-rat1}, we obtain
\begin{align*}
a \lambda_1^2 + b \bx_1 \cdot \bx_2 & = \frac{c \lambda_1^2}{2}\\
b \lambda_2^2 + a \bx_1 \cdot \bx_2 & = \frac{c \lambda_2^2}{2}.
\end{align*}
Notice that since \(\Lambda\) is not orthogonal, \(\bx_1\cdot \bx_2 \neq 0\) and 
$$\bx_1 \cdot \bx_2 = \frac{c-2a}{2b} \lambda_1^2 = \frac{c-2b}{2a} \lambda_2^2.$$
Now suppose that there are there are rational numbers \(p,q\) so that
$$p\lambda_1^2 = \bx_1\cdot\bx_2 = q\lambda_2^2.$$
Then, by \eqref{order-rat1}, there exist rational, and hence integer solutions $a,b,c$ to the linear system
\begin{equation}\label{order-mat-eq3}
\left\{ \begin{array}{ll}
    a\bx_1\cdot\bx_1 + b \bx_1\cdot \bx_2 + c \bx_1\cdot \bz = 0\\
    a\bx_1\cdot \bx_2 + b \bx_2\cdot\bx_2 + c \bx_2\cdot \bz = 0,
\end{array} \right.
\end{equation}
which factors as
\begin{equation}\label{order-mat-eq1}
    \begin{pmatrix}
    \bx_1^\top\\\bx_2^\top
    \end{pmatrix}
    \begin{pmatrix}
    \bx_1 &  \bx_2 & \bz
    \end{pmatrix}
    \begin{pmatrix}
    a\\b\\c
    \end{pmatrix}
    =
    \textbf{0}.
\end{equation}
Since the matrix $\begin{pmatrix} \bx_1^\top \\ \bx_2^\top \end{pmatrix}$ is of full rank, \((a,b,c)^\top\) solves \eqref{order-mat-eq1} if and only if it solves
\begin{equation}\label{order-mat-eq2}
    \begin{pmatrix}
    \bx_1 &  \bx_2 & \bz
    \end{pmatrix}
    \begin{pmatrix}
    a\\b\\c
    \end{pmatrix}
    =
    \textbf{0}.
\end{equation}
On the other hand, $\bz$ being an integer solution of \eqref{order-mat-eq2} is equivalent to \(\bz\) having finite order in \(\real^2/\Lambda\). In fact, the order of \(\bz\) in \(\real^2/\Lambda\) is \(\leq |c|\). By combining \eqref{order-mat-eq3} with \eqref{order-rat1}, we obtain the linear system
\begin{equation}
    \begin{pmatrix}
    2\lambda_1^2 & 2\bx_1\cdot\bx_2 & \lambda_1^2\\
    2\bx_1\cdot\bx_2 & 2\lambda_2^2 & \lambda_2^2
    \end{pmatrix}
    \begin{pmatrix}
    a\\b\\c
    \end{pmatrix}
    = \textbf{0}.
\end{equation}
If \(\lambda_1^2, \lambda_2^2, \bx_1 \cdot\bx_2 \in \zed\), then by Siegel's lemma (see, for instance, Lemma~4D in Chapter~1 of~\cite{schmidt1991}) there exists a nontrivial integer solution to this system with
$$\max \{ |a|,|b|,|c| \} \leq 3\sqrt{3} \left( \max \left\{ 2\lambda_1^2, 2 \lambda_2^2, |\bx_1 \cdot \bx_2| \right\} \right)^2 = 12 \sqrt{3}\ \lambda_2^4,$$
since \(\lambda_1^2 \leq | \bx_1\cdot \bx_2 | \leq \lambda_2^2\).
\endproof

\bigskip

\section{Equal covering extensions}
\label{cover}

In this section we investigate the covering radii of lattices in the plane, in particular proving Theorem~\ref{cover_ext} and its corollaries. Let $\Lambda \subset \real^2$ be a lattice of rank~$2$ with minimal basis $\bx,\bwy$ and angle $\theta \in [\pi/3,\pi/2]$ between these basis vectors. Then the successive minima of $\Lambda$ are
\begin{equation}
\label{min_basis}
0 < \lambda_1 = \|\bx\| \leq \lambda_2 = \|\bwy\|.
\end{equation}
See, for instance,~\cite{hex} for the details on the existence of such a minimal basis. We will use the following result in this section.

\begin{thm} [Theorem~3.2 of~\cite{cover-2}] \label{cover2} Consider the parallelogram generated by minimal basis vectors of $\Lambda$, as above. Then deep holes of $\Lambda$ in this parallelogram are the circumcenters of the two acute (right) triangles comprising this parallelogram and the covering radius of $\Lambda$ is the circumradius of the triangles.
\end{thm}

\noindent
We start with the following formula for the covering radius.

\begin{lem} \label{cover-1} The covering radius of $\Lambda$ is
\begin{equation}
\label{cv-eq}
\mu = \frac{\sqrt{\lambda_1^2 + \lambda_2^2 - 2\lambda_1\lambda_2\cos \theta}}{2 \sin \theta}.
\end{equation}
\end{lem}

\proof
The vectors $\bx,\bwy$ correspond to successive minima in~$\Lambda$, and hence form a reduced basis. Then Theorem~\ref{cover2} asserts that the covering radius of~$\Lambda$ is equal to the circumradius of the triangle with sides corresponding to the vectors $\bx$ and $\bwy$. The length of the third side of this triangle is
\begin{equation}
\label{c1}
\sqrt{\lambda_1^2 + \lambda_2^2 - 2 \lambda_1 \lambda_2 \cos \theta},
\end{equation}
and the area of this triangle is
\begin{equation}
\label{c2}
A = \frac{1}{2} \lambda_1 \lambda_2 \sin \theta.
\end{equation}
Now, the circumradius of a triangle with sides $a,b,c$ and area $A$ is given by the formula
\begin{equation}
\label{c3}
R = \frac{abc}{4A}.
\end{equation}
Putting together~\eqref{c1}, \eqref{c2} and \eqref{c3} produces~\eqref{cv-eq}.
\endproof

The similarity classes of WR lattices in the plane are parameterized by the angle $\theta \in [\pi/3,\pi/2]$, and each similarity class is represented by
$$L_{\theta} = \begin{pmatrix} 1 & \cos \theta \\ 0 & \sin \theta \end{pmatrix} \zed^2,$$
see~\cite{lf_pg_fl}, \cite{hex} for details. The following corollary follows immediately from Lemma~\ref{cover-1} by substituting $\lambda_1 = \lambda_2 = 1$ into~\eqref{cv-eq}.

\begin{cor} \label{cover_cor} The covering radius $\mu = \mu(\theta)$ of the lattices $L_{\theta}$ is a continuous function on the interval $[\pi/3,\pi/2]$, given by
$$\mu(\theta) = \frac{\sqrt{1 - \cos \theta}}{\sqrt{2} \sin \theta}.$$
The endpoints of the interval are represented by the hexagonal lattice and the square lattice $\zed^2$ with the covering radii  $1/\sqrt{3}$ and $1/\sqrt{2}$, respectively.
\end{cor}

We are now ready to prove Theorem~\ref{cover_ext}. We first want to build an extension $E_1 \subset \Lambda \subset \real^2$ with $\rk \Lambda = 2$ so that $\mu(\Lambda) = \mu(E_1)$. Our argument characterizes all possible such extensions, showing that they must be rectangular lattices, i.e. lattices containing an orthogonal basis.

\proof[Proof of Theorem~\ref{cover_ext}]
First notice that each $\Lambda(\alpha)$ as in~\eqref{cv_ext} is a rectangular  lattice, thus its successive minima are
$$\lambda_{1,2} = \sqrt{\alpha}, \sqrt{1 - \alpha},$$
i.e. norms of the orthogonal basis vectors given in~\eqref{cv_ext}. By Lemma~\ref{cover-1}, the covering of $\Lambda(\alpha)$ is 
$$\mu = \frac{\sqrt{\alpha + (1-\alpha) - 2 \sqrt{\alpha(1-\alpha)} \cos (\pi/2})}{2 \sin (\pi/2)} = \frac{1}{2}.$$

In the reverse direction, assume $\Lambda \subset \real^2$ is a full rank lattice so that $\be_1 \in \Lambda$ and $\mu(\Lambda) = 1/2$. The vector $\ba := \frac{1}{2} \be_1$ is a deep hole of $E_1$. First we show that $\ba$ is a deep hole of the lattice $\Lambda$ as well. Suppose not, then there exists a point $\bx = \begin{pmatrix} x_1 \\ x_2 \end{pmatrix} \in \Lambda$ such that
$$\| \bx - \ba \| < 1/2.$$
Then the vector $\bz = \be_1 - \bx = \begin{pmatrix} 1-x_1 \\ -x_2 \end{pmatrix}$ is also in $\Lambda$, and
$$\| \bz - \ba \| = \| \bx - \ba \| < 1/2.$$
Let $\Lambda' = \spn_{\zed} \{ \bx,\bz \} \subseteq \Lambda$, then $\mu(\Lambda') \geq \mu(\Lambda)$. The triangle with sides corresponding to the basis vectors $\bx, \bz$ of $\Lambda'$ is contained in the interior of the circle of radius $1/2$ with center at $\ba$, thus the circumradius $R$ of this triangle is $< 1/2$. On the other hand, by Theorem~\ref{cover2} the covering radius of~$\Lambda'$ is equal to the circumradius of the triangle with sides corresponding to the shortest basis vectors, which has to be $\leq R$. Hence $\mu(\Lambda) \leq \mu(\Lambda') < 1/2$, which is a contradiction, so $\ba$ is a deep hole of the lattice $\Lambda$. This means that there exists a basis $\bx,\bz \in \Lambda$ with $\|\bx\| = \lambda_1$, $\|\bz\| = \lambda_2$ so that the point $\ba$ is the center of the circle circumscribed around the triangle with sides $\bx,\bz$, meaning in particular that
\begin{equation}
\label{circle}
(1/2)^2 = \|\bx-\ba\|^2 = (x_1-1/2)^2 + x_2^2 = x_1^2 + x_2^2 - x_1 + (1/2)^2.
\end{equation}
Also, $2\ba = \be_1$ is a diagonal of the fundamental parallelogram of $\Lambda$ spanned by $\bx,\bz$, meaning that
$$\bx + \bz = \be_1.$$
Hence $\bz = \begin{pmatrix} 1- x_1 \\ -x_2 \end{pmatrix}$, and
$$\cos \theta = \frac{\bx \cdot \bz}{\|\bx\| \|\bz\|} = \frac{x_1 - x_1^2 - x_2^2}{\lambda_1 \lambda_2},$$
where $\theta$ is the angle between $\bx$ and $\bz$, which lies in the interval $[\pi/3,\pi/2]$. Hence $\cos \theta = 0$ by~\eqref{circle}. Letting $x_1 = \alpha$, we obtain
$$x_2 = \sqrt{\alpha - \alpha^2},$$
and~\eqref{cv_ext} follows by replacing $\bz$ with $-\bz$.
\medskip

Next, suppose that $L \subset \real^n$ be a lattice of rank~$1$ and let $\bu \in \real^n$ be such that $L = \zed \bu$, so the covering radius of $L$ is $\mu = \|\bu\|/2$. Let $\Lambda$ be a lattice of rank~$2$ in~$\real^n$ containing~$L$ and let $V = \spn_{\real} \Lambda$ be the $2$-dimensional subspace spanned by this lattice. Applying a suitable isometry of~$\real^n$, we can identify~$V$ with $C_2 := \{\bx \in \real^n : x_i = 0\ \forall\ 2 < i \leq n\}$. In fact, we can choose such an isometry $\tau$ so that $\bu$ maps to $\beta \be_1$ for $\beta = \|\bu\|$. Then $\Lambda' = \frac{1}{\beta} \tau(\Lambda)$ is a lattice isometric to $\frac{1}{\beta} \Lambda$ in $\real^n$, and $\Lambda'$ contains $\be_1$. Identifying $C_2$ with $\real^2$ we see that Theorem~\ref{cover_ext} implies that $\Lambda'$ is an equal covering extension of $\zed \be_1$ in $\real^2$ if and only if it is of the form~\eqref{cv_ext}. Finally, notice that
$$\det(L) = \sqrt{\det (\bu^{\top} \bu)} = \|\bu\| = \beta.$$
This completes the proof of the theorem.
\endproof

\begin{rem} \label{rem_WR}  An immediate implication of Theorem~\ref{cover_ext} is that the only well-rounded equal covering extension of $E_1$ is 
\begin{equation}
\label{Lambda12}
\Lambda(1/2) = \begin{pmatrix} 1/2 & -1/2 \\ 1/2 & 1/2 \end{pmatrix} \zed^2,
\end{equation}
which is a square lattice in the plane containing $\zed^2$ as a sublattice of index~$2$. More generally, a rank-two equal covering extension $\Lambda \subset \real^n$ of a rank-one lattice $L \subset \Lambda$ is well-rounded if and only if it is isometric to $\det(L) \Lambda(1/2)$. Further, the set of all similarity classes of planar lattices is parameterized by
$$F = \{ (a,b) \in \real^2 : 0 \leq a \leq 1/2, a^2+b^2 \geq 1 \},$$
see Figure~\ref{domain}. The set of semi-stable classes in~$\real^2$ contains the WR classes: from~\eqref{ss_det} its follows that a lattice $\Lambda$ in~$\real^2$ is semi-stable if and only if $\lambda_1(\Lambda) \geq \det (\Lambda)^{1/2}$ (see~\cite{lf_pg_fl} for more details). Thus the only semi-stable equal covering extensions are also those similar to~$\Lambda(1/2)$ as in~\eqref{Lambda12}, i.e. similar to~$\zed^2$ as demonstrated in Figure~\ref{domain}.
\begin{figure}
\centering
\includegraphics[scale=0.4]{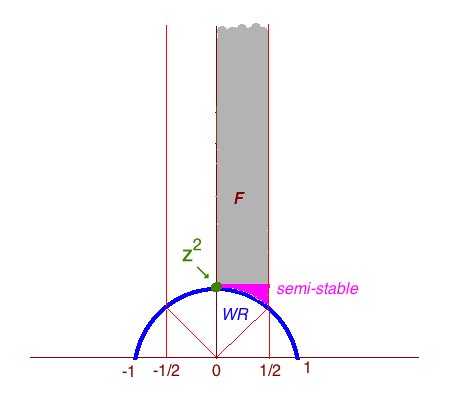}
\caption{Similarity classes of planar lattices with $\zed^2$ representing the only equal covering extension class that is WR and semi-stable.}\label{domain}
\end{figure}
\end{rem}

Let $D$ be a squarefree integer and $K = \que(\sqrt{D})$ be a quadratic number field with embeddings $\sigma_1,\sigma_2 : K \to \cee$. Let $\Sigma_K : K \to \real^2$ be the Minkowski embedding of $K$, defined for every $x \in K$ as
$$\Sigma_K(x) = \begin{pmatrix} \sigma_1(x) \\ \sigma_2(x) \end{pmatrix},$$
if $D > 0$ and
$$\Sigma_K(x) = \begin{pmatrix} \Re(\sigma_1(x)) \\ \Im(\sigma_1(x)) \end{pmatrix},$$
if $D < 0$. We write $\Omega_K$ for the planar lattice $\Sigma_K(\O_K)$, where $\O_K$ is the ring of integers of the number field $K$.

\begin{cor} \label{nf_cover-1} Assume $D \not\equiv 1\ (\md 4)$, then $\Omega_K$ is an equal covering extension of the rank-one lattice~$\zed \Sigma_K(1+\sqrt{D})$.
\end{cor}

\proof
If $D \not\equiv 1\ (\md 4)$ then $\O_K = \zed[\sqrt{D}]$, and so
$$\Omega_K = \begin{pmatrix} 1 & \sqrt{D} \\ 1 & -\sqrt{D} \end{pmatrix} \zed^2 \text{ if } D > 0,\ \Omega_K = \begin{pmatrix} 1 & 0 \\ 0 & \sqrt{|D|} \end{pmatrix} \zed^2 \text{ if } D < 0.$$
In either case, the lattice $\Omega_K$ is rectangular. If $D > 0$, then $\lambda_1 = \sqrt{2}$, $\lambda_2 = \sqrt{2D}$, and so Lemma~\ref{cover-1} implies that $\mu(\Omega_K) = \sqrt{D+1}/\sqrt{2}$, while
$$\Sigma_K(1+\sqrt{D}) = \begin{pmatrix} 1 + \sqrt{D} \\ 1 - \sqrt{D} \end{pmatrix},$$
and so $\mu \left(\zed \Sigma_K(1+\sqrt{D}) \right) = \sqrt{D+1}/\sqrt{2}$. If $D < 0$, then $\lambda_1 = 1$, $\lambda_2 = \sqrt{|D|}$, and so Lemma~\ref{cover-1} implies that $\mu(\Omega_K) = \sqrt{|D|+1}/2$, while
$$\Sigma_K(1+\sqrt{D}) = \begin{pmatrix} 1 \\ \sqrt{|D|} \end{pmatrix},$$
and so $\mu \left(\zed \Sigma_K(1+\sqrt{D}) \right) = \sqrt{|D|+1}/2$.
\endproof

\begin{cor} \label{nf_cover-2} Assume $D \equiv 1\ (\md 4)$, then $\Omega_K$ is not an equal covering extension of any rank-one lattice.
\end{cor}

\proof
If $D \equiv 1\ (\md 4)$ then $\O_K = \zed \left[ \frac{1+\sqrt{D}}{2} \right]$, and so
$$\Omega_K = \begin{pmatrix} 1 & \frac{1+\sqrt{D}}{2} \\ 1 & \frac{1-\sqrt{D}}{2} \end{pmatrix} \zed^2 \text{ if } D > 0,\ \Omega_K = \begin{pmatrix} 1 & \frac{1}{2} \\ 0 & \frac{\sqrt{|D|}}{2} \end{pmatrix} \zed^2 \text{ if } D < 0.$$
In both cases it is not difficult to check that $\Omega_K$ does not have an orthogonal basis, and hence cannot be similar to a lattice of the form $\Lambda(\alpha)$ as in~\eqref{cv_ext}. The conclusion follows from Theorem~\ref{cover_ext}.
\endproof

Finally, we discuss a construction of orthogonal equal covering extensions in any dimension.

\proof[Proof of Theorem~\ref{mu_orth}]
We will argue by induction on $k \geq 1$. Theorem~\ref{cover_ext} establishes the base of induction, so let \(k \geq 2\). Let \(\{\bx_1,\ldots, \bx_k\}\subset \real^n\) be an orthogonal basis for \(\Lambda_k\) and let \(\be_{k+1} \in \real^n\) be a vector orthogonal to \(\spn_{\real} \Lambda_k\). Let 
\[
  P_k = \left\{\sum_{i=1}^k a_i\bx_i: a_i \in \{0,1\}, 1 \leq i\leq k\right\}
\]
be the set of vertices of the fundamental parallelepiped spanned by \(\bx_1,\ldots, \bx_k\). The circumcenter of this orthogonal parallelepiped is the point \(\bz \in \real^n\) which is equidistant from the points of \(P_k\) by Lemma~\ref{simplex}, and hence is a deep hole of $\Lambda_k$. Let 
\[
  B_k = \left\{ \bwy \in \spn_\real\{\bx_1, \ldots, \bx_k\}: \|\bwy-\bz\| = \mu(\Lambda_k) \right\}.
\]
Let \(\Lambda_{k-1} = \spn_\zed \{\bx_1,\ldots, \bx_{k-1}\}\) and let 
\[
  P_{k-1} = \left\{\sum_{i=1}^{k-1} a_i\bx_i: a_i \in \{0,1\}, 1 \leq i\leq {k-1}\right\}.
\]
Now define
\[
  B_{k-1} = B_k \cap \spn_\real \Lambda_{k-1}
\]
By construction, \(P_{k-1}\subset B_{k-1}\), while \(B_{k-1}\) is a sphere in a \((k-1)\)-dimensional subspace and the points of \(P_{k-1}\) are elements of an orthogonal lattice in that subspace. Let \(\bz'\) be the orthogonal projection of \(\bz\) onto \(\spn_\real (\Lambda_{k-1})\).
Since \(\{\bx_1,\ldots,\bx_k\}\) is an orthogonal set, \(\bz' = \bz - \text{proj}_{\bx_k}(\bz)\).
Moreover, \(\bz\) is equidistant from \(\{\bx_1,\ldots, \bx_k\}\), so by Lemma~\ref{simplex}
\[
  \begin{pmatrix}
  \bx_1^\top\\\vdots\\\bx_k^\top
  \end{pmatrix}
  \bz
  =\frac{1}{2} 
  \begin{pmatrix}
  \|\bx_1\|\\\vdots\\\|\bx_k\|
  \end{pmatrix},
\]
and therefore
\[
  \begin{pmatrix}
  \bx_1^\top\\\vdots\\\bx_{k-1}^\top
  \end{pmatrix}
  \bz'
  =\frac{1}{2} 
  \begin{pmatrix}
  \|\bx_1\|\\\vdots\\\|\bx_{k-1}\|
  \end{pmatrix}.
\]
Then \(\bz'\) is equidistant from \(\{\bx_1,\ldots,\bx_{k-1}\}\) and \(\Lambda_{k-1}\) is an orthogonal lattice contained in the $k$-dimensional subspace \(V = \spn_{\real} \{\Lambda_{k-1}, \be_{k+1} \}\). By the induction hypothesis, there exists a rank \(k\) orthogonal lattice \(\Lambda_k' \subset V\) so that \(\Lambda_{k-1}\subset \Lambda_k'\) and \(\bz'\) is a deep hole of $\Lambda_k'$. Let \(\textbf{y}_1,\ldots,\textbf{y}_k\) be an orthogonal basis for \(\Lambda_k'\) so that \(\textbf{y}_1,\ldots,\textbf{y}_k\) are equidistant from \(\bz'\). Since \(\bz = \bz'- \text{proj}_{\bx_k}(\bz)\), and \(\bx_k\) is orthogonal to \(V\),  \(\textbf{y}_1,\ldots,\textbf{y}_k\) are also equidistant from \(\bz\). Let \(\Lambda_{k+1} = \spn_\zed \{\textbf{y}_1,\ldots,\textbf{y}_k, \bx_k\}\). Then \(\Lambda_{k+1}\) is an orthogonal lattice that contains \(\Lambda_k\) and 
\[
  \begin{pmatrix}
  \textbf{y}_1^\top\\\vdots\\\textbf{y}_{k-1}^\top\\\bx_k^\top
  \end{pmatrix}
  \bz
  =\frac{1}{2} 
  \begin{pmatrix}
  \|\textbf{y}_1\|\\\vdots\\\|\textbf{y}_{k-1}\|\\\|\bx_k\|
  \end{pmatrix}.
\]
Thus \(\bz\) is equidistant from $\bo,\bwy_1,\dots,\bwy_k$ by Lemma~\ref{simplex}, and hence a deep hole of~\(\Lambda_k\).
\endproof
\medskip

\noindent
{\bf Acknowledgement.} We would like to thank the two anonymous referees for a very thorough reading of our paper and many helpful suggestions that improved the quality of exposition.
\bigskip

\bibliographystyle{plain}  

\end{document}